\documentclass[twocolumn]{article}
\usepackage[utf8]{inputenc}
\usepackage[left=2cm, right=2cm, top=2cm]{geometry} 
\usepackage[small]{titlesec}
\usepackage{etoolbox}
\usepackage{amsthm}
\makeatletter
\patchcmd{\ttlh@hang}{\parindent\z@}{\parindent\z@\leavevmode}{}{}
\patchcmd{\ttlh@hang}{\noindent}{}{}{}
\makeatother

	\usepackage{mathptmx} 
	\usepackage{times} 
	\usepackage{amsmath} 
	\usepackage[cmintegrals]{newtxmath}
	\usepackage{algorithm}
	\usepackage{algpseudocode}
	\usepackage{color}
	\usepackage{bm}
	\usepackage{makecell}
	\usepackage{subfig, float}
	\usepackage{tikz,pgfplots}
	\usepackage{epstopdf}
	\usepackage{natbib}
	\usepackage{mathtools}
	\usepackage[colorlinks=true,linkcolor=black,anchorcolor=black,citecolor=black,filecolor=black,menucolor=black,runcolor=black,urlcolor=black]{hyperref}
	\hypersetup{breaklinks=true}
	\usepackage{accents}

	\usetikzlibrary{external,arrows,shapes,positioning,plotmarks, arrows.meta, automata}
	\usepgfplotslibrary{external} 
	
\newlength\figureheight 
\newlength\figurewidth  

\DeclareMathOperator*{\subjecttoo}{subject\:to}

\DeclareMathAlphabet{\mathcalOld}{OMS}{cmsy}{m}{n}
\DeclareMathOperator*{\argmin}{arg\,min}

\DeclareMathOperator*{\minimizer}{minimize}
\newcommand{\tfinal}{t_f}
\newcommand{\tinitial}{t_0}
\newcommand{\xinitial}{\bm{x}_0}
\newcommand{\xfinal}{\bm{x}_f}

\newcommand{\alphabm}{\bm{\alpha}}

\newcommand{\nbm}{\bm{n}}
\newcommand{\xbm}{\bm{x}}
\newcommand{\ubm}{\bm{u}}
\newcommand{\taubm}{\bm{\tau}}
\newcommand{\sship}{\mathcalOld{S}_b}

\newcommand{\vship}{\mathcalOld{V}_b}

\newcommand{\sspatial}{\mathcalOld{S}_{\text{env}}}

\newcommand{\ubar}[1]{\underaccent{\bar}{#1}}

\title{\LARGE \bf
A COLREGs-Compliant Motion Planner for Autonomous Maneuvering of Marine Vessels in Complex Environments \footnote{This work was partially supported by FFI/VINNOVA and the Wallenberg Artificial Intelligence, Autonomous Systems and Software Program (WASP) funded by Knut and Alice Wallenberg Foundation; and Saab AB.}
}

\date{}
\author{Kristoffer Bergman \thanks{Kristoffer Bergman and Daniel Axehill are with the Division of Automatic Control,
		Link{\"o}ping University, Link\"oping, Sweden
		(e-mail: {\tt\small \{kristoffer.bergman, daniel.axehill\}@liu.se})} \and Oskar Ljungqvist  \thanks{Oskar Ljungqvist is at Volvo Trucks, Gothenburg, Sweden, e-mail: \texttt{oskar.ljungqvist@volvo.com.}}  \and Jonas Linder \thanks{Jonas Linder is at ABB AB, Corporate Research, V\"aster{\aa}s, Sweden (e-mail:
		\tt\small jonas.x.linder@se.abb.com)} \and Daniel Axehill \footnotemark[2]  
}

\begin{document}
	
	\algblock{ParFor}{EndParFor}
	\algnewcommand\algorithmicparfor{\textbf{parfor}}
	\algnewcommand\algorithmicpardo{\textbf{do}}
	\algnewcommand\algorithmicendparfor{\textbf{end\ parfor}}
	\algrenewtext{ParFor}[1]{\algorithmicparfor\ #1\ \algorithmicpardo}
	\algrenewtext{EndParFor}{\algorithmicendparfor}

	\maketitle
	\thispagestyle{empty}
	\pagestyle{empty}

		\textbf{\textit{Abstract ---}}\textbf{An enabling technology for future sea transports is safe and energy-efficient autonomous maritime navigation in narrow environments with other marine vessels present. This requires that the algorithm controlling the ship is able to account for the vessel's dynamics, and obeys the rules specified in the international regulations for preventing collision at sea (COLREGs). To account for these properties, this work proposes a two-step optimization-based and COLREGs-compliant motion planner tailored for large vessels. In the first step, a lattice-based motion planner is used to compute a suboptimal trajectory based on a library of precomputed motion primitives. To comply with the rules specified in COLREGs, the lattice-based planner is augmented with discrete states that represent what type of COLREGs situations that are active with respect to other nearby vessels. The trajectory computed by the lattice-based planner is then used as warm-start for the second receding-horizon step based on direct optimal control techniques. The aim of the second step is to solve the problem to local optimality with the discrete COLREGs states remained fixed. The applicability of the proposed motion planner is demonstrated in simulations, where it computes energy-efficient and COLREGs-compliant trajectories. The results also show that the motion planner is able to prevent complex collision situations from occurring.    }
\section{Introduction} \label{sec:intro}
The interest in autonomous ship technology has over the past decade rapidly increased. The main driving forces in the development of future sea transports are to improve efficiency, e.g., to minimize the environmental footprint by reducing fuel consumption, and to increase safety by reducing the number of accidents caused by human errors~\citep{riviera2020industry, imo2020imo}. An enabling technology to achieve these goals is autonomous maritime navigation in narrow environments such as archipelagos, when there are other manned and unmanned marine vessels present. As a consequence, it is important that the algorithm controlling the ship is able to avoid collisions with both static and dynamic obstacles\footnote{Further on in this paper, the term ``ship'' will be used to denote the ship that is being controlled, and ``dynamic obstacles'' to denote other marine vessels covered by COLREGs.}. In maritime navigation, the rules specified by the international regulations for preventing collisions at sea (COLREGs) must be followed when there is a risk of collision between vessels~\citep{cockcroft2003guide}. Hence, it is required that algorithms in the control system comply with these rules.   

In the literature, there exist a number of works that have focused on collision avoidance algorithms which comply with COLREGs. In~\citep{kuwata2013safe}, a motion planner is developed based on so-called velocity obstacles. The same concept is used in~\citep{kufoalor2018proactive}, where it is extended to include cooperative behavior. Another COLREGs-compliant motion planner is presented in \citep{hagen2018mpc}, where a receding-horizon approach is used based on input sampling to a closed-loop controller. In \citep{eriksen2020hybrid}, a hierarchical approach for collision avoidance is proposed, where a trajectory is computed in a first step using optimal control without considering dynamic obstacles. In a second step, a receding horizon controller is used to track the nominal trajectory while complying with COLREGs. Common for \citep{kuwata2013safe,kufoalor2018proactive,hagen2018mpc,eriksen2020hybrid} are that they all use some sort of nominal reference that is computed without considering dynamic obstacles.

A probabilistic approach is to use the closed-loop variant of the rapidly-exploring random tree (RRT) algorithm for COLREGs-compliant motion planning~\citep{chiang2018colreg}, where the input to an underlying controller is sampled to generate feasible trajectory candidates. This planner is shown to efficiently compute COLREGs-compliant trajectories in several scenarios but is unable to optimize the trajectories with respect to a performance measure. 

A method based on a COLREGs-compliant repairing A$^\star$ graph search is proposed in \citep{campbell2014automatic}. The algorithm first computes a path while neglecting the system dynamics, and  then smoothen the path using splines. This method works well for agile unmanned surface vehicles (USVs). However, to be able to compute safe and energy-efficient trajectories for large vessels, a motion planner that is able to account for the vessel's dynamics is required. There are many papers that propose optimization-based motion planners for computing safe and energy-efficient trajectories for vessels, e.g.,~\citep{martinsen2019autonomous, martinsen2020optimization, bitar2020two}. However, neither of these works consider COLREGs. 

This work proposes a two-step optimization-based motion planner for marine vessels that complies with COLREGs. It is based on the two-step algorithm presented in~\citep{bergman2020marine} and is here extended to account for dynamic obstacles while complying with COLREGs. In the first step, a lattice-based motion planner~\citep{pivtoraiko2009differentially} is used to compute a suboptimal trajectory to a discretized version of the motion planning problem. By employing classical graph-search
algorithms, the lattice-based planner concatenates motion
segments from a finite library of optimized \textit{motion primitives}~\citep{bergman2019improved}. This step is here augmented with discrete states that represent what types of COLREGs situations that are active for the surrounding dynamic obstacles in each time step. 
To improve the trajectory computed by the lattice-based planner, a second receding-horizon step based on direct optimal control techniques is added that solves the problem to local optimality while keeping the discrete COLREGs states fixed. Together, the two steps constitute the proposed COLREGs-compliant motion planning algorithm. In an online setting, an update in the predicted trajectories of the dynamic obstacles is allowed to trigger re-planning.  

The proposed framework allows to use a high-fidelity \textit{maneuvering} model of the vessel~\citep{fossen2011handbook}, making it possible to compute safe and energy-efficient trajectories in narrow environments. Compared to other approaches for COLREGs-compliant motion planning that modifies a trajectory based on a nominal reference, the proposed method uses an economic approach and computes the nominal trajectory based only on the initial and terminal states and the objective function. This will potentially prevent complex collision situations from occurring at all.

The remainder of the paper is organized as follows. In Section~\ref{sec:model} the dynamical model of the ship is briefly described. In Section~\ref{sec:colregs}, the COLREGs-compliant motion planning problem is presented, and the proposed framework for solving the problem is given in Section~\ref{sec:motionplanning}. A simulation study is presented in Section~\ref{sec:simulations} and the paper is concluded in Section~\ref{sec:conclusions} with a summary of the presented work as well as some suggestions for potential future research. 

\section{Ship modeling }\label{sec:model}
\begin{figure}[!t]
	\centering
	\setlength\figureheight{0.25\textwidth}
	\setlength\figurewidth{0.4\textwidth}
%
%
\definecolor{mycolor1}{rgb}{0.89020,0.10196,0.10980}%
\definecolor{mycolor2}{rgb}{0.2000,0.6275,0.1725}%

\begin{tikzpicture}

\begin{axis}[%
width=\figurewidth,
height=\figureheight,
at={(0\figurewidth,0\figureheight)},
scale only axis,
xmin=-63.2602939004499,
xmax=39.8511971894683,
ymin=-12.8581762917281,
ymax=12.8786287775784,
axis background/.style={fill=white},
axis x line*=bottom,
axis y line*=left,
axis line style={draw=none},
tick style={draw=none},
ytick=\empty,
xtick=\empty,
axis equal,
]
\addplot [color=mycolor1, line width=1.4pt, forget plot]
table[row sep=crcr]{%
	37.2193889287494	-7.31755117800078\\
	39.5511971894683	-19.7755985947342\\
	28.1856742996503	-25.3849804361991\\
	-41.2602939004499	9.33800366385103\\
	-32.2265792713508	27.4054329220493\\
	37.21938892874940	-7.31755117800078\\
};

\addplot[area legend, draw=black, fill=white!91!black, forget plot]
table[row sep=crcr] {%
	x	y\\
	32.5236461759999	-5.52869679600094\\
	34.1261166902488	-6.52560441982046\\
	35.5233590554552	-7.80529729666178\\
	36.6729187693914	-9.32889260601816\\
	37.5398670418544	-11.050096675391\\
	38.0978620882828	-12.9166115901039\\
	38.32994951155	-14.8717242389754\\
	38.2290774527414	-16.8560295135041\\
	37.7983108583204	-18.8092353020304\\
	37.7263747566166	-18.9531075054381\\
	36.4222700824481	-20.4696442540907\\
	34.8953490981103	-21.7409250658547\\
	33.1920065250528	-22.7283227165639\\
	31.3639975654254	-23.401835628249\\
	29.4668653464494	-23.7409994519022\\
	27.5582532706025	-23.7355088663671\\
	25.6961535500056	-23.3855307003068\\
	23.9371451424007	-22.7016988631993\\
	23.9371451424007	-22.7016988631993\\
	-40.5894735072	9.56161046160101\\
	-32.0029724736008	26.7346125287994\\
	32.5236461759999	-5.52869679600095\\
	32.5236461759999	-5.52869679600094\\
}--cycle;
\addplot [color=black, draw=none, mark=x, mark options={solid, black}, forget plot]
table[row sep=crcr]{%
	-28.6216701119973	14.3108350559987\\
};

\node[] (source) at (axis cs:-29.6,16.3){};
\node (destination) at (axis cs:-22.3,0.15){};
\draw[mycolor2,->, ultra thick, >=latex](source)--(destination);
\node[] (source) at (axis cs:-30.6,15.3){};
\node (destination) at (axis cs:-0,0){};
\draw[dotted, ultra thick, >=latex](source)--(destination);
 \draw[  thick, ->] (axis cs:-19,9) arc [radius=12,start angle=-26.56,end angle=-60];

\node[] (source) at (axis cs:-64.1,-30){};
\node (destination) at (axis cs:-30,-30){};
\draw[blue,->, ultra thick, >=latex](source)--(destination);
\node[] (source) at (axis cs:-62,-32){};
\node (destination) at (axis cs:-62,0){};
\draw[blue,->, ultra thick, >=latex](source)--(destination);
\node[blue] at (axis cs: -30,-27) { $y_E$};
\node[blue] at (axis cs: -58,-2) { $x_E$};
\node[black] at (axis cs: -15,2) { $\alpha_j$};
\node[mycolor2] at (axis cs: -50,18) { $\bm{\tau}_j(\alpha_j, n_{j}, \bm{\nu})$};
\node[] (source) at (axis cs:0,-1){};
\node (destination) at (axis cs:0,20){};
\draw[dotted, ultra thick, >=latex](source)--(destination);
\node[] (source2) at (axis cs:-26.83/11,13.42/11){};
\node (destination2) at (axis cs:26.83,-13.42){};
\draw[blue, ->, ultra thick, >=latex](source2)--(destination2);
 \draw[  thick, ->] (axis cs:0,3) arc [radius=3,start angle=90,end angle=-26.56];
\node[] (source) at (axis cs:13.42/11,26.83/11){};
\node (destination) at (axis cs:-13.42,-26.83){};
\draw[blue,->, ultra thick, >=latex](source)--(destination);
\node[] at (axis cs: 5,2.2) { $\psi$};
\node[blue] at (axis cs: 27,-10) { $x_B$};
\node[blue] at (axis cs: -16,-20) { $y_B$};
\node[] at (axis cs: 27,-16.) { $u$};
\node[] at (axis cs: -14,-27) { $v$};

\node[] (source) at (axis cs:-5,26){};
\node (destination) at (axis cs:-5,12){};
\draw[->,thick, >=latex](source)--(destination);
\node[] at (axis cs: -5,27.0) { $\sship(\xbm)$};
\end{axis}
\end{tikzpicture}%
 \vspace*{-0.1em}
	\caption{ \label{fig:ship} \small Definition of Earth-fixed ($x_E,y_E$) and body-fixed ($x_B, y_B$) coordinate systems, body-fixed velocities ($u,v$) and the convex polytope $\sship(\xbm)$ representing the bounding region of the ship's body at $\xbm$. Furthermore, the force vector of the $j^{\mathrm{th}}$ thruster is illustrated.   }
\end{figure}
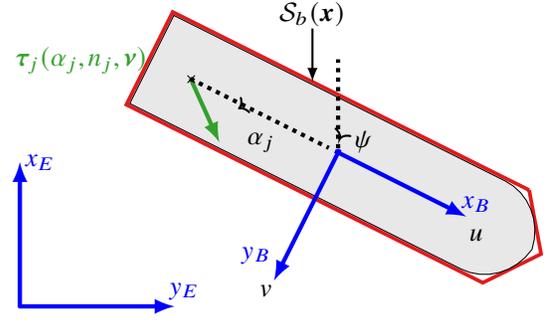

This section contains a brief description of the dynamic ship model used for motion planning. 
The ship is assumed to move on the ocean surface and it is thus sufficient to consider the horizontal 3 degrees of freedom motion. The motion of the ship is described using two coordinate systems: an Earth-fixed (inertial) system and a body-fixed system that is located at the center of mass of the ship. The Earth-fixed generalized position is given by $\bm{\eta} = \begin{bmatrix}
x  &y  &\psi\end{bmatrix}^\intercal \in  \mathbf{SE}(2)$~\citep{lavalle2006planning} and the body-fixed generalized velocity vector is represented by $\bm{\nu} = \begin{bmatrix}u  &v  &r\end{bmatrix}^\intercal \in \mathbf{R}^3$, see Fig.~\ref{fig:ship}. The generalized velocity and position is related through the kinematics
\begin{equation}\label{eq:kinematics}
\dot{\bm{\eta}} = \bm{R}(\psi)\bm{\nu}.
\end{equation}
where $\bm{R}(\psi)$ is the rotation matrix. The motion induced by forces acting on the ship is derived using rigid-body mechanics and theory of hydrodynamics, see \citep{fossen2011handbook} for details. The model is given by
\begin{equation} \label{eq:dynamics}
\bm{M}\dot{\bm{\nu}} + \bm{C}(\bm{\nu}) \bm{\nu} + \bm{D}(\bm{\nu}) \bm{\nu} = \bm{\tau},
\end{equation}
where $\bm{M}$ is the total inertia matrix, $\bm{C}(\bm{\nu}) \bm{\nu}$ corresponds to Coriolis and centripetal forces, $\bm{D}(\bm{\nu})\bm{\nu}$ describes the damping and $\bm{\tau}$ represents the  forces that are acting on the ship. Simplified models for the thrusters are used where it is assumed that forces induced by the propellers and rudders can be separated. The total forces can be written as 
\begin{equation} \label{eq:thrusters}
\taubm(\bm{\alpha}, \bm{n}, \bm{\nu}) = \sum_{j=1}^{N_t} \bm{T}_j \taubm_j(\alpha_j, n_j, \bm{\nu}),
\end{equation}
where $N_t$ is the number of thrusters, $\bm{T}_j$ describes the geometry of the thruster configuration for thruster $j$ and the vectors \mbox{$\alphabm = \begin{bmatrix}
	\alpha_1  &\ldots  &\alpha_{N_t}\end{bmatrix}^\intercal$} and \mbox{$\nbm = \begin{bmatrix}
	n_1  &\ldots  &n_{N_t}\end{bmatrix}^\intercal$} represent the thruster angles and velocities. For more details regarding thruster and ship modeling, the reader is referred to, e.g.,~\citep{bergman2020marine, fossen2011handbook,perez2006ship}.

To account for dynamics in the thrusters, the derivatives of the thruster angles and velocities are considered as control inputs to the system. A compact notation of the state vector is thus defined as $\bm{x} = \begin{bmatrix}\bm{\eta}^\intercal \; \; \bm{\nu}^\intercal \; \; \alphabm^\intercal \; \; \nbm^\intercal \end{bmatrix}^\intercal$
and the control-input vector is \mbox{$\bm{u} = \begin{bmatrix}\dot{\alphabm}^\intercal &\dot{\nbm}^\intercal \end{bmatrix}^\intercal$}. 
Then by combining \eqref{eq:kinematics}-\eqref{eq:thrusters}, the ship dynamics can be written as
\begin{align} \label{eq:dynamics-summarized}
\dot{\bm{x}} = \bm{f}({\bm{x}},{\bm{u}}).
\end{align}
The feasible set for the state space is given by 
\begin{align}
\mathcalOld{X} = \{ \xbm \; | \; ||\begin{bmatrix}
u &v\end{bmatrix}|| \leq v_{\mathrm{max}}, \; |n_j| \leq \bar{n}, \quad j = 1, \ldots, N_t \},
\end{align} where $v_{\mathrm{max}}$ is the maximum allowed magnitude of the ship's velocity and $\bar{n}$ is the maximum speed for the thrusters. The control inputs are restricted to
\begin{align}
\mathcalOld{U} = \{\ubm \; | \; |\dot{\alpha}_j| \leq \bar{\dot{\alpha}}, \quad |\dot{n}_j| \leq \bar{\dot{n}}, \quad j = 1 ,\ldots, N_t \},
\end{align}  
where $\bar{\dot{\alpha}}$ and $\bar{\dot{n}}$ are the maximum angular rate and acceleration of the thrusters, respectively. 

The ship is assumed to operate in an environment with both static and dynamic obstacles. The bounding region of the ship's body at state $\xbm$ is defined by $ \sship(\xbm)$, see Fig.~\ref{fig:ship}. The region occupied by static obstacles is given by $\mathcalOld{X}^s_{\mathrm{obst}}$, and the time-varying region for dynamic obstacles is $\mathcalOld{X}^d_{\mathrm{obst}}(t)$. This region can, e.g., be described as the union of convex polytopes which depend on the states of the dynamic obstacles at time~$t$. Thus, it is possible to pose the complete obstacle region as the union of static and dynamic obstacles
\begin{equation}  \label{eq:obstacle_constraints}
\mathcalOld{X}_{\mathrm{obst}}(t) = \mathcalOld{X}^s_{\mathrm{obst}} \cup \mathcalOld{X}^d_{\mathrm{obst}}(t).
\end{equation}
The free space where the vessel is not in collision with any obstacle is defined as $\mathcalOld X_{\mathrm{free}}(t) = \mathcalOld X \setminus \mathcalOld X_{\mathrm{obst}}(t)$.
Since this set is the complement set of $\mathcalOld X_{\mathrm{obst}}(t)$, it is in general not~convex.  

\section{COLREGs-compliant motion planning}\label{sec:colregs}
This section contains the description of a mathematical representation of the COLREGs regulations, which is used to formulate the COLREGs-compliant motion planning problem. 

\subsection{COLREGs representation}
\begin{figure}[t!]
	\centering
	\setlength\figureheight{0.18\textwidth}
	\setlength\figurewidth{0.46\textwidth}
	\input{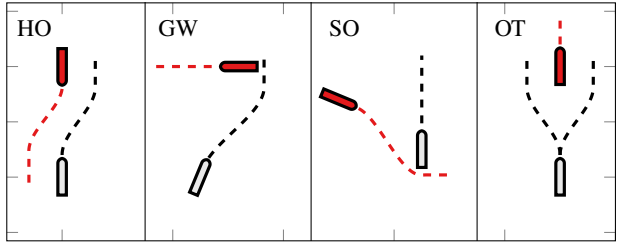}
	\caption{ \small \label{fig:illustrative_example} The main COLREGs scenarios together with the conceptually correct vessel behaviors. The ship is grey and the obstacle vessel is red. From left: Head on (HO), give way (GW), stand on (SO) and overtaking (OT).    }
\end{figure}
The regulations covered by COLREGs are divided into several parts and there are in total 41 rules. For the purpose of planning trajectories for collision avoidance, part B is most relevant. The rules covered in this work are Rule 8 and Rule 13-17, see~\citep{cockcroft2003guide} for details. The corresponding COLREGs situations that are handled are 
\begin{itemize}
	\item \textbf{Safe (SF)}. The vessels are far away or heading away from each other.
	\item \textbf{Head on (HO)}. The vessels are meeting on a near reciprocal course with a risk of collision.
	\item \textbf{Give way (GW)}. Crossing where the dynamic obstacle approaches from the starboard side. 
	\item \textbf{Stand on (SO)}. Crossing where the dynamic obstacle approaches from the port side, or when being overtaken.
	\item \textbf{Overtaking (OT)}. The ship is about to overtake the dynamic obstacle.
	\item \textbf{Emergency (EM)}. When the dynamic obstacle is not complying with COLREGs.
\end{itemize} 
Some of these situations are illustrated in Fig.~\ref{fig:illustrative_example}. A discrete state $q^i(t) \in \mathcalOld{Q}$ is used to represent which of these situations that is active for the ship with respect to dynamic obstacle $i$, where $\mathcalOld{Q} = \{\text{SF}, \text{HO}, \text{OT}, \text{GW}, \text{SO}, \text{EM}\}$. 

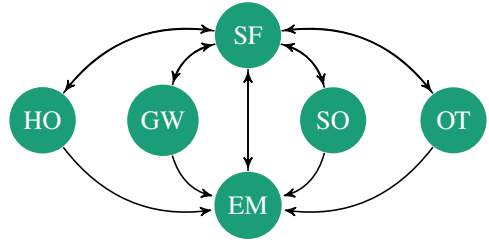
\begin{figure}[t!]
	\centering
	\definecolor{mycolor1}{RGB}{27,158,119}%
\definecolor{mycolor2}{RGB}{117,112,179}%
\begin{tikzpicture}[->,>=stealth',shorten >=1pt,auto,node distance=1.6cm,
                    semithick]
  \tikzstyle{every state}=[fill=mycolor1,draw=none,text=white]

  \node[state] 		   (A)                    {SF};
  \node[state]         (B) [below left of=A] {GW};
  \node[state]         (C) [left of=B]       {HO};
  \node[state]         (D) [below right of=A] {SO};
  \node[state]         (E) [right of=D] {OT};
  \node[state]         (F) [below left of=D] {EM};

  \path (A) edge [bend right]             node {} (B)
            edge [bend right]             node {} (C)
            edge [bend left]             node {} (D)
            edge [bend left]             node {} (E)
            edge                         node {} (F)
        (B) edge [bend right] 			 node {} (F)
        	edge [bend left] 			 node {} (A)
        (C) edge [bend right] 			 node {} (F)
        	edge [bend left] 			 node {} (A)
        (D) edge [bend left] 			 node {} (F)
        	edge [bend right] 			 node {} (A)
        (E) edge [bend left] 			 node {} (F)
        	edge [bend right] 			 node {} (A)
        (F) edge 			 			 node {} (A);
\end{tikzpicture}
	\caption{ \small \label{fig:transitions} Allowable switches between COLREGs states.    }
\end{figure}

Each dynamic obstacle is modeled using a finite-state machine. The allowable switches between states in $\mathcalOld{Q}$ are illustrated in Fig.~\ref{fig:transitions}. Common choices for switching criterias are to use measures based on closest point of approach (CPA)~\citep{kuwata2013safe,eriksen2020hybrid,kufoalor2018proactive}. The time to CPA is calculated by assuming that the ship and dynamic obstacle keep their current velocity $\bm{v} = \begin{bmatrix}u &v\end{bmatrix}^\intercal$ constant. Define $\bm{p} = \begin{bmatrix} x &y\end{bmatrix}^\intercal$ and $\bm{v}$ as the position and velocity of the ship, and $\bm{p}_o$ and $\bm{v}_o$ as the position and velocity of the dynamic obstacle. Then, it is possible to calculate the time to CPA according to 
\begin{equation}
t_{\text{CPA}} = \begin{cases}
0 & \lVert \bm{v} - \bm{v}_o \rVert \leq \epsilon, \\
\frac{(\bm{p}_o - \bm{p})\cdot(\bm{v} - \bm{v}_o)}{\lVert \bm{v} - \bm{v}_o \rVert^2 } & \text{otherwise}.
\end{cases}
\end{equation}
This expression is derived by solving the equation \mbox{$\frac{\delta}{\delta t} \lVert (\bm{p} + t\bm{v}) - (\bm{p}_o + t\bm{v}_o) \rVert = 0$} with respect to $t$ \citep{kufoalor2018proactive}. The distance between the ship and the dynamic obstacle at CPA~is 
\begin{equation}
d_{\text{CPA}} = \lVert (\bm{p} + t_{\text{CPA}}\bm{v}) - (\bm{p}_o + t_{\text{CPA}}\bm{v}_o) \rVert.
\end{equation}

Similar to what is done in~\citep{eriksen2020hybrid}, a switch from state SF to either HO, GW, SO or OT (see Fig~\ref{fig:transitions}) occurs when
\begin{equation} \label{eq:switch}
\left(d_{\text{CPA}} \leq \ubar{d} \right) \; \land \; \left( \ubar{t}_{\text{low}} \leq t_{\text{CPA}} \leq \ubar{t}_{\text{high}} \right),  
\end{equation} 
where $\ubar{d}, \; \ubar{t}_{\text{low}}, \; \ubar{t}_{\text{high}}$ are threshold parameters which can be state dependent. If $\eqref{eq:switch}$ is true, a switch to a new COLREGs state occurs. As in~\citep{tam2010collision}, the COLREGs state is selected based on the relative position and heading of the dynamic obstacle with respect to the ship.  A switch back to state SF occurs when 
\begin{equation} \label{eq:switch-back}
\left( d_{\text{CPA}} \leq \bar{d} \right) \; \lor \;  \lnot \left( \bar{t}_{\text{low}} \leq t_{\text{CPA}} \leq \bar{t}_{\text{high}} \right),  
\end{equation}
where $\ubar{d} < \bar{d}$, $\bar{t}_{\text{low}} < \ubar{t}_{\text{low}}$ and $\ubar{t}_{\text{high}} < \bar{t}_{\text{high}}$ are threshold parameters that implement hysteresis and thus avoid fast switches between states.
Note that it is not possible to switch between any of the situations described in Fig.~\ref{fig:illustrative_example} without first switching to SF. This is according to what is described in the regulations~\citep{cockcroft2003guide}.

A switch to EM only occurs during re-planning if any of the surrounding dynamic obstacles do not comply with COLREGs, or if a new dynamic obstacle is detected that is dangerously close to the ship. As in~\citep{eriksen2020hybrid}, a critical distance $d_{\text{crit}}$ is defined to handle switches to EM. The time to the critical point $t_{\text{crit}}$ is found by solving 
\begin{equation}
\lVert (\bm{p} + t_{\text{crit}}\bm{v}) - (\bm{p}_o + t_{\text{crit}}\bm{v}_o) \rVert = d_{\text{crit}}.
\end{equation}    
If $t_{\text{crit}} < \ubar{t}_{\text{crit}}$, where $\ubar{t}_{\text{crit}}$ is a user-defined threshold parameter, the state is switched to EM. Note in Fig.~\ref{fig:transitions} that it is possible to switch to EM from any COLREGs state, as opposed to what is implemented in~\citep{eriksen2020hybrid},  

For notational brevity, the switching between states are written as 
\begin{equation}
q^i(t) = \sigma(\xbm(t), \; \xbm^i_o(t), \; q^i(t^-)),
\end{equation}  
where $\sigma$ is the switching function\footnote{The switching function is allowed to depend on the type of dynamic obstacle.}, while $\xbm(t)$ and $\xbm^i_o(t)$ are the states of the ship and dynamic obstacle $i$, respectively. The notation $q^i(t^-)$ is used for left-hand limit of $q^i$ at time~$t$.

One commonly used approach for COLREGs-compliant collision avoidance is to assign a spatial region for each obstacle that should not be entered by the ship. The definitions of these spatial regions are not covered by the rules in COLREGs, and there exist several different implementation approaches to obtain the desired behavior. Some examples are to use cost functions~\citep{eriksen2020hybrid, hagen2018mpc}, forbidden regions~\citep{tam2010collision, campbell2014automatic}, or velocity obstacles~\citep{kuwata2013safe, kufoalor2018proactive}. 
This work uses cost functions, but the method presented in Section~\ref{sec:motionplanning} can be extended to account for any of the above mentioned approaches. 

The cost functions are all based on the relative position of the ship with respect to dynamic obstacles. The relative position between a set $\mathcalOld{S}_1 \subset \mathbf{R}^2$ and a position $ \bm{y} \in \mathbf{R}^2$ is defined as
\begin{equation}
\text{rel\_pos}(\mathcalOld{S}_1, \bm{y} )= \underset{\Delta \bm{p}}{\argmin}  \{ \lVert \Delta \bm{p} \rVert \; : \; \bm{x} - \Delta \bm{p} = \bm{y}, \; \bm{x} \in \mathcalOld{S}_1  \}. 
\end{equation}
The COLREGs related cost function is defined as
\begin{equation} \label{eq:colregs-ca}
\ell_{\text{CA}}(\xbm, \bm{x}^i_{o}, q^i) = \begin{cases}
k_{q^i} + \ell_{q^i}(\Delta \bm{p}^i_o, \xbm_o^i) & \hspace{-0.1em} \Delta \bm{p}^i_o \in \mathcalOld{S}_{\text{CA}}(q^i, \xbm_o^i) \\
k_{q^i} & \hspace{-0.1em} \text{otherwise}
\end{cases}
\end{equation}
where $\Delta \bm{p}^i_o = \text{rel\_pos}(\sship(\xbm), \bm{p}^i_o)$ is the relative position of the ship with respect to obstacle $i$. Moreover, $\ell_{q^i} > 0$ is the cost function related to COLREGs state $q^i$ at state $\xbm_o^i$ and \mbox{$\mathcalOld{S}_{\text{CA}}(q^i, \xbm_o^i) \subset \mathbf{R}^2$} is the set of relative positions for which the cost function is active. Finally, $k_{q^i}$ is added such that it is possible to penalize the time spent in the state $q^i$. 
\subsection{Problem formulation}
The COLREGs-compliant motion planning problem (CCMP) is now defined as follows: compute a feasible and collision-free state and control input trajectory $(\bm x(t),\bm u(t))$, $t\in[0,T_f]$, that moves the ship from its initial state \mbox{$\bm{x}(0)=\xinitial \in \mathcalOld X_{\mathrm{free}}$} to a desired goal state ${\bm{x}}(T_f)=\xfinal\in \mathcalOld X_{\mathrm{free}}$, while minimizing an objective functional $J$. The cost function used to define $J$ is separated into three parts. The first part is $\ell(\xbm,\ubm)$, which can, e.g., be related to time and energy consumption. The second part is given by $\ell_s(\xbm)$ and
penalizes the distance to static obstacles. Finally, the third part is defined in \eqref{eq:colregs-ca} and ensures that the COLREGs rules are not violated for the $N_o$ dynamic obstacles. 
The CCMP can be posed as the following continuous-time hybrid OCP:
\begin{alignat}{2} 	\label{eq:ocp}
&\minimizer_{\ubm(\cdot), \tfinal} 
&& \hspace{0em}\int_{\tinitial}^{\tfinal} \left[ \ell(\xbm, \ubm) + \ell_s(\xbm)  ) + \sum^{N_o}_{i=1} \ell_{\mathrm{CA}}(\xbm, \bm{x}^i_o, q^i)) \right] \mathrm{d}t  \nonumber \\
& \subjecttoo \hspace{1em} 
&& \xbm(\tinitial) = \xinitial, \quad \xbm(\tfinal)= \xfinal, \quad q^i(\tinitial) = q^i(\tinitial^-)\nonumber \\ 
& \text{} && \dot{\xbm}(t) = \bm{f}(\xbm(t), \ubm(t)), \nonumber \\
& \text{} && q^i(t) = \sigma(\xbm(t), \xbm^i_o(t), q^i(t^-) ), \quad i = 1 \ldots N_o \\
& \text{} && \xbm(t) \in \mathcalOld{X}_{\mathrm{free}}(t), \; \ubm(t) \in \mathcalOld{U}.\nonumber 
\end{alignat} 
The upcoming section describes how the proposed framework for solving \eqref{eq:ocp} is implemented by using a combination of a lattice-based motion planner and an optimization-based improvement step.

\section{Optimization-based motion planner} \label{sec:motionplanning}
In this section, the two steps of the proposed motion planner are presented in detail, and an algorithm for COLREGs-compliant motion planning is outlined.

\subsection{Lattice-based motion planner} \label{sec:lbmp}
A lattice-based motion planner computes a suboptimal solution to the continuous-time hybrid OCP in~\eqref{eq:ocp}. It is done by discretizing and transforming \eqref{eq:ocp} to a graph-search problem which is solved using suitable techniques. The transformation is done by restricting the control-inputs to a finite subset of the available actions. In this work, the so-called state-lattice methodology is used~\citep{pivtoraiko2009differentially}, where the search space is discretized into $\mathcalOld{X}_d$ and the set of available actions is represented using a motion-primitive set $\mathcalOld{P}$. 

A motion primitive $\bm{m}$ is defined as a feasible state and control-input trajectory that moves the ship from an initial state $\xbm_0 \in \mathcalOld{X}_d$ to a terminal state $\xbm_f \in \mathcalOld{X}_d$ in time $T$. By discretizing the time uniformly at $t_j = 0 + j\Delta_{t_m}, \; j = 0, \ldots N_m$, a sampled version of a motion primitive at these time instances is defined as
\begin{equation} \label{eq:motionPrimitive}
\bm{m} =  \{ \xbm_0, \; \ubm_0, \; \xbm_1, \ubm_1, \ldots, \xbm_{N_m} \}.	\end{equation}
To compute the motion primitives $\bm{m} \in \mathcalOld{P}$, the optimization-based framework in~\citep{bergman2019improved} is used. Since the positions of the obstacles are unknown offline, the motion primitives are computed by solving OCPs where only a discretized version of the position invariant cost
\begin{equation} \label{eq:pos-inv-cost-primitive}
\ell_{p,m}(\bm{m})= \Delta_{t_m}\sum_{j=0}^{N_m-1} \ell(\xbm_j, \ubm_j)
\end{equation} 
is minimized.  More details on the motion primitive computation for the ship is found in~\citep{bergman2020marine}.
Together, the motion primitives define the motion primitive set $\mathcalOld{P}$. Since the ship is position invariant, each motion primitive $\bm{m} \in \mathcalOld{P}$ can be translated and reused from each position in $\mathcalOld{X}_d$. The state transition function is defined as $\xbm_j = \bm{f}_m(\bar{\xbm}, \bm{m}, t_j)$ which outputs the state after $t_j$ seconds when motion primitive $\bm{m}$ is applied from $\bar{\xbm}$. 

The obstacle dependent cost functions are evaluated online during the graph search. When a motion primitive $\bm{m}$ is applied from state $\bar{\xbm}$, the cost related to static obstacles is
\begin{equation}
\ell_{s,m}(\bar{\xbm}, \bm{m}) = \Delta_{t_m} \sum_{j=0}^{N_m-1} \ell_s( \xbm_j),
\end{equation}
where $\xbm_j =\bm{f}_m(\bar{\xbm}, \bm{m}, t_j)$. If the location of the static obstacles are known offline, it is possible to precompute an approximation of this function in terms of a static cost map. For the dynamic obstacles, it is assumed that a prediction of the state of the obstacle $\bm{x}_o^i(\cdot)$ is given, such that the COLREGs state trajectory $\bm{q}^i = \{q^i_{0}, \ldots q^i_{N_m-1} \}$ can be computed for a motion primitive $\bm{m}$ applied at state $\bar{\xbm}$. This is used to define the cost related to dynamic obstacles for $\bm{m}$ applied at $\bar{\xbm}$:
\begin{equation} \label{eq:dyn-cost-primitive}
\small
\ell_{\text{CA},m}\left(\bar{\xbm}, \bm{m}, \bm{x}^i_o(\cdot), \bm{q}^i\right) = \Delta_{t_m} \hspace{-0.5em}\sum^{N_m-1}_{j=0}  \ell_{\text{CA}}\left(\xbm_j, \bm{x}^i_o(\bar{t}_0+t_j), q^i_{j}\right),
\end{equation} 
where  $\xbm_j =\bm{f}_m(\bar{\xbm}, \bm{m}, t_j)$ and $\bar{t}_0$ is the time when the motion primitive is applied. One approach to approximate this function online is to compute the relative position to $\xbm^i_o(\bar{t}_0+t_j)$ for a finite number of vertices in the ownship's body $\sship(\bm{x}_j)$. 

By summarizing the terms in \eqref{eq:pos-inv-cost-primitive}--\eqref{eq:dyn-cost-primitive}, the total stage cost when a motion primitive is applied is 
\begin{equation} \label{eq:lat-cost}
\begin{aligned}
&\ell_m\left(\bar{\xbm}, \bm{m}, \{\bm{x}_o^i(\cdot) \}^{N_o}_{i=1}, \{\bm{q}^i \}_{i=1}^{N_o} \right) = \\ 
&\ell_{p,m}(\bm{m}) + \ell_{s,m}(\bar{\xbm}, \bm{m}) + \sum_{i=1}^{N_o} \ell_{\text{CA},m}\left(\bar{\xbm}, \bm{m}, \bm{x}^i_o(\cdot), \bm{q}^i\right).
\end{aligned}
\end{equation}
Now, the resulting graph-search problem can be posed as a discrete OCP in the form:
\begin{subequations}
	\label{eq:lattice}
	\begin{alignat}{2}
	&\minimizer_{\{ \bm{m}_k \}_{k=1}^{M},\hspace{0.5ex}M}
	&&  \sum_{k=1}^{M} \ell_m\left(\bar{\xbm}_k, \bm{m}_k, \{\bm{x}_o^i(\cdot) \}^{N_o}_{i=1}, \{\bar{\bm{q}}_k^i \}_{i=1}^{N_o} \right)  \nonumber\\
	&\hspace{0.1em} \subjecttoo \hspace{1em} 
	&&  \bar{\xbm}_1 = \xinitial, \quad \bar{\xbm}_{M+1} = \xfinal, \quad \bar{\bm{q}}^i_0 = q_i(\tinitial^-)  \\ \label{eq:primStateTrans}
	& \text{} &&  \bar{\xbm}_{k+1} = \bm{f}_{m}(\bar{\xbm}_{k}, \bm{m}_{k}, t_{N_{m_k}}), \\ \label{eq:hehe}
	& \text{} && \bar{\bm{q}}^i_{k} = \bm{\sigma}_{m}(\bar{\xbm}_k, \bm{m}_{k}, \xbm^i_{o}(\cdot), \bar{\bm{q}}^i_{k-1}) \\ \label{eq:primInSet}
	& \text{} &&  \bm{m}_{k}  \in \mathcalOld{P}(\xbm_k), \\ \label{eq:trajInFreeSpace}
	& \text{} &&  c(\bar{\xbm}_k, \bm{m}_k) \notin \mathcalOld{X}_{\mathrm{obst}},
	\end{alignat} 
\end{subequations}
where the decision variables are the motion-primitive sequence $\{ \bm{m}_k \}_{k=1}^{M}$ and its length $M$. The function in \eqref{eq:hehe} outputs the COLREGs state trajectory for dynamic obstacle~$i$ along $\bm{m}_k$ applied from $\bar{\xbm}_k$, and \eqref{eq:primInSet} restricts the choice of $\bm{m}_k$ to the set of applicable motion primitives $\mathcalOld{P}(\bar{\xbm}_k)$ from state $\bar{\xbm} _k$. Finally, the constraint in \eqref{eq:trajInFreeSpace} ensures that the ship does not collide with obstacles when $\bm{m}_k$ is applied from $\bar{\xbm}_k$. 

Since the problem in \eqref{eq:lattice} is in the form of a graph-search problem, standard  graph-search algorithms can be used to solve it~\citep{pivtoraiko2009differentially}. For a faster search, a precomputed free-space heuristic look-up table (HLUT)~\citep{knepper2006high} is used as heuristic function. The output is a state and control-input trajectory ($\bar{\xbm}(\cdot), \bar{\ubm}(\cdot) )$ which is defined by the motion primitive sequence. Additional outputs are the COLREGs state trajectories for all obstacles. For dynamic obstacle $i$, it is defined by 
\begin{equation} \label{eq:latcolreg}
\bar{q}^i(t) = q^i_k \quad \ubar{t}_k \leq t < \bar{t}_k, \quad k = 1, \ldots, N_k
\end{equation}
where $\ubar{t}_k$ and $\bar{t}_k$ with $\ubar{t}_k = \bar{t}_{k-1}$ define the time interval where COLREGs state $q^i_k$ is active, and $N_k$ 
is the number of COLREGs state switches.

\subsection{Receding-horizon improvement} \label{sec:rhi}
The aim of the receding-horizon improvement step is to reduce the suboptimality introduced from solving \eqref{eq:lattice}. This is achieved by performing local improvements to the trajectories from the lattice-based planner $(\bar{\xbm}(\cdot), \bar{\ubm}(\cdot))$ in a receding-horizon fashion over a user-defined planning horizon~$T$. This step is based on the work in~\citep{bergman2020marine}, and is here extended to account for dynamic obstacles and COLREGs. 


By using direct methods for optimal control, the improvement OCP is constructed by discretizing the continuous-time motion planning problem~\eqref{eq:ocp}, e.g., using direct multiple-shooting with $N+1$ discretization points~\citep{diehl2006fast}, where the time between two points is $\Delta_t = T/N$. 
To avoid collisions with both static and dynamic obstacles, a local spatial constraint is computed at each discrete point $j$, $j = 1, \ldots, N-1$. This constraint is represented by a convex polytope $\sspatial^j$ which is defined as
\begin{equation} \label{eq:sspatial}
\sspatial^j = \{ \bm{p} \in \mathbf{R}^2 \; | \; \bm{A}_j \bm{p} \leq \bm{b}_j \},  
\end{equation}
where $\bm{A}_j = [\bm{a}_{j,1} \; \ldots \; \bm{a}_{j,K_j} ]^\intercal \in \mathbf{R}^{K_j \times 2}$, $\bm{a}_{j,l} \in \mathbf{R}^2, \; ||\bm{a}_{j,l}||_2 = 1,$ for $l = 1 ,\ldots ,K_j$, $\bm{b}_j \in \mathbf{R}^{K_j}$ and $K_j$ is the number of half spaces that defines $\sspatial^j$. These polytopes are computed using Algorithm~1 in \citep{bergman2020marine}.
Since both $\sship(\xbm_j)$ and $\sspatial^j$ are convex, a sufficient condition to ensure collision avoidance for the ship's body is that all vertices $\vship$ of $\sship(\xbm_j)$ lie within $\sspatial^j$~\citep{martinsen2019autonomous}. This condition can be described by the following constraint:
\begin{equation} \label{eq:ocp_coll}
\bm{A}_j\bm{c}_{\text{rot}}(\xbm_j, \bm{p}_l) \leq \bm{b}_j, \quad \forall \bm{p}_l \in \vship,
\end{equation}
where 
\begin{align}
\bm{c}_{\text{rot}}(\xbm_j, \bm{p}_l) =\begin{bmatrix}
\cos \psi_j & -\sin \psi_j \\ \sin \psi_j & \cos \psi_j
\end{bmatrix}\bm{p}_l + \begin{bmatrix}
x_j \\ y_j
\end{bmatrix},
\end{align} and $\bm{A}_j$ and $\bm{b}_j$ are the half-space representation of $\sspatial^j$ defined in~\eqref{eq:sspatial}. Furthermore, we want to express the position-dependent cost related to the distance to obstacles also in the improvement step. Here, this is done by modifying \eqref{eq:ocp_coll} to:
\begin{equation} \label{eq:ocp_coll_epsi}
\begin{aligned}
\bm{A}_j\bm{c}_{\text{rot}}(\xbm_j, \bm{p}_l) &\leq \bm{b}_j - \bm{1}_j(d_{\mathrm{safe}}-\epsilon_{d,j}), \quad \forall \bm{p}_l \in \vship,
\end{aligned}
\end{equation}
where a variable $\epsilon_{d,j}$ that satisfies $0 \leq \epsilon_{d,j} \leq d_{\mathrm{safe}}$ is added, which represents the distance to the boundaries of $\sspatial^j$ up to a user-defined safety distance $d_{\mathrm{safe}}$. With this modification, the objective function in the improvement step is 
\begin{equation} \label{eq:imp-cost}
\begin{aligned}
&\ell_{\text{imp}}\left(\xbm, \ubm, \epsilon_{d}, \{\bm{x}^i_o,  q^i\}^{N_o}_{i=1}\right) = \\
&\ell(\xbm, \ubm) + k_d \epsilon_{d}^2 + \sum^{N_o}_{i=1} \ell_{\mathrm{CA}}(\xbm, \bm{x}^i_o. q^i).	
\end{aligned}
\end{equation}
From the current state $\xbm_{\text{cur}}$ at time $t_{\text{cur}}$, the receding-horizon improvement OCP is
\begin{alignat}{2} 	\label{eq:improvement}
&\minimizer_{\hspace{0.7em}\{\bm{u}_j, \epsilon_{d,j}\}^{N-1}_{j=0} } 
&& \hspace{0.7em}  \Delta_t \sum_{j=0}^{N-1} \left[ \ell_{\text{imp}}(\xbm_j, \ubm_j, \epsilon_{d,j}, \{\bm{x}^i_o(t_j),  \bar{q}^i(t_j)\}^{N_o}_{i=1})   \right]   \nonumber \\
& \hspace{0.7em}\subjecttoo 
&& \hspace{0.7em} \xbm_0 = \xbm_{\text{cur}}, \quad \xbm_N= \bar{\xbm}(t_{\text{cur}}+T), \quad t_j = t_{\text{cur}} + j \Delta_t, \nonumber \\ 
& \text{} && \hspace{0.7em} \xbm_{j+1} = \bm{f}_d(\xbm_j, \ubm_j, \Delta_t),  \\
& \text{} && \hspace{0.7em} 	\bm{A}_j\bm{c}_{\text{rot}}(\xbm_j, \bm{p}_l)  \leq \bm{b}_j - \bm{1}_j(d_{\mathrm{safe}} - \epsilon_{d,j}), \; \forall \bm{p}_l \in \vship \nonumber\\
& \text{} && \hspace{0.7em} \xbm_j \in \mathcalOld{X}, \; \ubm_j \in \mathcalOld{U}, \quad 0 \leq \epsilon_{d,j} \leq d_{\mathrm{safe}}.\nonumber 
\end{alignat} 
Here, $\bm{f}_d(\xbm, \ubm, \Delta_t)$ is a discretized version of \eqref{eq:dynamics-summarized} over the time interval $\Delta_t$. Furthermore, the terminal state $\xbm_N$ is constrained such that the solution is connected to the initial trajectory $\bar{\xbm}(\cdot)$ at the end of the planning horizon, which is done for stability guarantees~\citep{bergman2020optimization}. Compared to the initial problem formulation~\eqref{eq:ocp}, the combinatorial aspect of selecting the COLREGs state trajectories $\bar{q}^i(\cdot)$ is in the improvement step pre-determined by the lattice-based planner as defined in \eqref{eq:latcolreg}. Furthermore, the second combinatorial aspect of selecting on which side to pass obstacles is implicitly encoded in the warm-start trajectory from the lattice planner. Hence, the hybrid OCP in~\eqref{eq:ocp} has been transformed into a standard discrete-time OCP, which can be solved using direct OCP techniques such as the ones presented in~\citep{diehl2006fast}. 

\subsection{Proposed CCMP algorithm} \label{sec:pccmpa}

\begin{algorithm}[t]
	\caption{Proposed CCMP algorithm}
	\label{alg:ccmp}
	\begin{algorithmic}[1]
		\State \textbf{input}: $\xinitial$, $\xfinal$, $\mathcalOld{P}$, $\mathcalOld{X}^s_{\text{obst}}$,  $\{\xbm_o^i(\cdot) \}^{N_o}_{i=1}$, $T$, $\Delta_t$
		\State $\left( \bar{\xbm}(\cdot), \bar{\ubm}(\cdot), \{\bar{q}^i(\cdot) \}_{i=1}^{N_o}  \right)  \leftarrow $ Solve \eqref{eq:lattice} with $\mathcalOld{P}$ \label{line:lattice}
		\State $t_\text{cur} \leftarrow \tinitial + \Delta_t$
		\While {$\bar{\xbm}(t_\text{cur} + T) \neq \xfinal$} \label{line:while}
		\State $\{\xbm_o^i(\cdot) \}^{N_o}_{i=1} \leftarrow$ \texttt{obstacles\_prediction}() \label{line:obst}
		\State $\bar{q}^i(t) \leftarrow \sigma\left(\bar{\xbm}(t), \bar{\xbm}_o^i(t), \bar{q}^i(t^-) \right), \;  t_\text{cur} \leq t \leq \tfinal  $  \label{line:colreg}
		\If {\textbf{any} $\bar{q}^i(t_\text{cur})$ = EM }
		\State \texttt{emergency\_maneuver}($\bar{\xbm}(t_\text{cur})$)
		\ElsIf {\textbf{any} $\bar{\xbm}(t) \in \mathcalOld{X}^d_{\text{obst}}(t), \; t \geq t_{\text{cur}} + T$ } \label{line:collision}
		\State $\left( \bar{\xbm}(\cdot), \bar{\ubm}(\cdot), \{\bar{q}^i(\cdot) \}_{i=1}^{N_o}  \right)   \leftarrow $ Solve \eqref{eq:lattice} from $ \bar{\xbm}(t_\text{cur})$
		\Else
		\State $\{\bm{A}_j, \bm{b}_j \}_{j=1}^{N-1} \leftarrow $ \texttt{prepare\_improvement}($\bar{\xbm}(\cdot)$)
		\State $\left( \bar{\xbm}(\cdot), \bar{\ubm}(\cdot) \right) \leftarrow $ Solve \eqref{eq:improvement} from $\xbm_\text{cur} = \bar{\xbm}(t_\text{cur})$ \label{line:improvement}
		\State re\_plan $\leftarrow$ \texttt{check\_trajectory}($\bar{\xbm}(\cdot)$)
		\If { re\_plan \textbf{is} True}
		\State \hspace{-1em}$\left( \bar{\xbm}(\cdot), \bar{\ubm}(\cdot), \{\bar{q}^i(\cdot) \}_{i=1}^{N_o}  \right)   \leftarrow $ Solve \eqref{eq:lattice} from $ \bar{\xbm}(t_\text{cur})$
		\EndIf
		\EndIf
		\State Send trajectory to controller : 
		
		\hspace{1.5ex} \texttt{send\_reference}($\left( \bar{\xbm}(\cdot), \bar{\ubm}(\cdot) \right)$ \label{line:traj-track}
		\State $t_\text{cur} \leftarrow t_\text{cur}  + \Delta_t$
		\EndWhile
	\end{algorithmic}
\end{algorithm}
In an online setting, an update in the predicted trajectories of the dynamic obstacles might require a re-planning. To handle such scenarios, a CCMP algorithm is proposed in Algorithm~\ref{alg:ccmp} that is based on the two steps described in \mbox{Section~\ref{sec:lbmp}} and Section~\ref{sec:rhi}. 

In the first step on line~\ref{line:lattice}, the lattice-based motion planner is used to compute the initial state and control-input trajectories ($\bar{\xbm}(\cdot), \bar{\ubm}(\cdot)$) and COLREGs state trajectories $\{\bar{q}^i(\cdot) \}_{i=1}^{N_o} $ by solving \eqref{eq:lattice}. Then, $(\bar{\xbm}(\cdot), \bar{\ubm}(\cdot))$ is improved in a receding-horizon fashion until the terminal state is reached.  
Due to uncertainties in the prediction of the dynamic obstacle trajectories $\{\xbm_o^i(\cdot) \}_{i=0}^{N_o}$, and that new dynamic obstacles might appear during execution, the predictions are updated at the beginning of each improvement iteration (line~\ref{line:obst}). Then, the associated COLREGs state trajectories are updated on line~\ref{line:colreg} according to the new obstacle predictions. If any of the COLREGs state trajectories are in state EM, an emergency maneuver from the current state is computed. There are several ways to compute such a maneuver, e.g., to use the lattice-based motion planner for computing a safe stop maneuver based on the current situation.

If the updated dynamic obstacle predictions cause a collision beyond the current improvement planning horizon $t_{\text{cur}} + T$ (line~\ref{line:collision}), a new trajectory is computed by the lattice-based planner.
Otherwise, a preparation step is performed where the initial and terminal constraints are specified, the objective function is parameterized based on $\bar{q}^i(t)$ and the collision avoidance polytopes $\bm{A}_j, \bm{b}_j$ are computed. After the preparation step, \eqref{eq:improvement} is solved from $\xbm_{\text{cur}}$ on line~\ref{line:improvement}. If the trajectory from solving \eqref{eq:improvement} triggers re-planning, e.g., if it is close to infeasibility or the objective function value has increased due to that the dynamic obstacle predictions have changed, the lattice-based motion planner computes a new trajectory from $\xinitial = \bar{\xbm}(t_{\text{cur}})$. Finally, the state and control-input trajectories are sent on line~\ref{line:traj-track} to a trajectory-tracking controller such as~\citep{barlund2020nonlinear} to stabilize the ship and suppress various disturbances.  

\section{Simulation study} \label{sec:simulations}
\begin{figure}[t!]
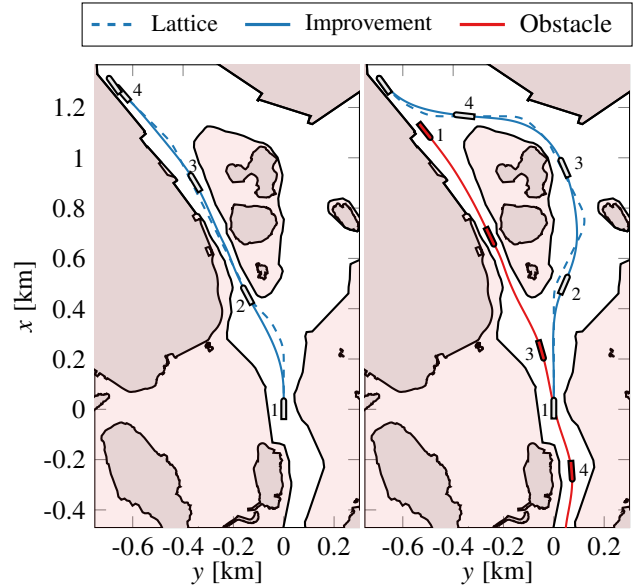

	\centering
	\hspace*{1.5em}\definecolor{mycolor1}{rgb}{0.6510,0.8078,0.8902}%
\definecolor{mycolor2}{rgb}{0.1216,0.4706,0.7059}%
\definecolor{mycolor3}{rgb}{0.6980,0.8745,0.5412}%
\definecolor{mycolor4}{rgb}{0.2000,0.6275,0.1725}%
\definecolor{mycolor5}{rgb}{0.9843,0.6039,0.6000}%
\definecolor{mycolor6}{rgb}{0.8902,0.1020,0.1098}%
\definecolor{mycolor7}{rgb}{ 0.9922,0.7490,0.4353}%
\definecolor{mycolor8}{rgb}{1.0000,0.4980,0}%
\definecolor{mycolor9}{rgb}{ 0.7922,0.6980,0.8392}%
\definecolor{mycolor10}{rgb}{0.4157,0.2392,0.6039}%
\begin{tikzpicture}

\begin{axis}[%
hide axis,
width=0,
height=0,
at={(0,0)},
scale only axis,
xmin=0,
xmax=1,
ymin=0,
ymax=1,
axis background/.style={fill=white},
xmajorgrids,
ymajorgrids,
legend style={legend cell align=left,column sep=3pt, align=center},
legend columns=-1,
]

\addplot [color=mycolor2, dashed, line width=1.2pt]
table[row sep=crcr]{%
	1	0\\
};
\addlegendentry{ {\small Lattice } }

\addplot [color=mycolor2, line width=1.2pt]
table[row sep=crcr]{%
	1	0.0\\
};
\addlegendentry{ {\small Improvement  } }

\addplot [color=mycolor6, line width=1.2pt]
table[row sep=crcr]{%
	1	0\\
};
\addlegendentry{ Obstacle }

\end{axis}
\end{tikzpicture}%
 \\[-0.15cm]
	\hspace*{-2em}  
	\subfloat[No dynamic obstacle (P1)\label{subfig-1:dummy}]{%
		\setlength\figureheight{0.35\textwidth}
		\setlength\figurewidth{0.2\textwidth}
		\input{sim_avoid_nom.tex}\vspace*{-2em}
	}
	\subfloat[Dynamic obstacle (P2)\label{subfig-2:dummy}]{%
		\hspace{-1.5em} 
		\setlength\figureheight{0.35\textwidth}
		\setlength\figurewidth{0.2\textwidth}
		\input{sim_avoid.tex} 
	}
	\caption{ \small \label{fig:sim_avoid} Two similar docking scenarios where the white area represents navigable water. The vessels are plotted at given timestamps (1-4) for interpretabilty. In (a), the motion planner finds a trajectory through the narrow passage. In (b), a dynamic obstacle occupies the narrow passage, which makes the proposed motion planner to find a trajectory on the other side of the island such that it completely avoids a complex COLREGs situation.   }
\end{figure}

In this section, the proposed framework is evaluated in a simulation study by solving motion planning problems in a model of the Helsinki archipelago\footnote{Video available at \texttt{\url{https://youtu.be/yPEFJRywWHQ}}}. The geometry of the environment is represented using polygons that have been calculated using data from Open Street Map and the Finnish Transport and Communications Agency (Traficom). These polygons have been expanded using depth data such that only navigable water is represented, see, e.g., Fig.~\ref{fig:sim_avoid}. The cost map used in the lattice-based planner is computed with resolution 1 m. The motion planner is implemented in C++, where the OCPs are solved using CasADi~\citep{andersson2018casadi} with IPOPT~\citep{wachter2006implementation} as NLP solver, and discretized using multiple shooting. 
From the results in~\citep{bergman2020marine}, it was found that a planning horizon of \mbox{$T = 150$ s} with \mbox{$\Delta_t = 2$ s} are suitable values to use for the supply vessel with respect to the trade-off between computation time and optimality. These values are used in all scenarios presented in this section when solving the improvement OCP in~\eqref{eq:improvement}, which means that the improvement step runs at a frequency of 0.5 Hz.

\begin{figure}[t]
	\centering
	\setlength\figureheight{0.18\textwidth}
	\setlength\figurewidth{0.46\textwidth}
	\input{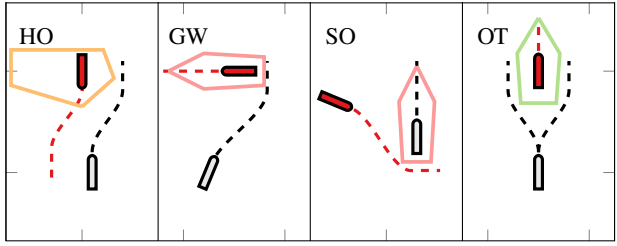}
	\caption{ \small \label{fig:colregs_borders} Different polygons used to represent the restricted regions associated with the COLREGs scenarios shown in Fig.~\ref{fig:illustrative_example}. }
\end{figure}
The simulations are based on a model of a small ship with two thrusters. The model of the ship is derived from the \emph{supply vessel} model in the MSS hydro toolbox~\citep{perez2006overview}. The parameters for the two thrusters are:
\begin{equation}
\bar{n} = 2 \; \mathrm{RPS} \quad \bar{\dot{n}} = 0.08 \; \mathrm{RPS/s}, \quad \bar{\dot{\alpha}} = 7.2 ^\circ\mathrm{/s},	
\vspace*{-0.5em}
\end{equation}
The lattice-based planner is discretized with a position resolution of $r_p = 5$ m, and the maximum velocity is set to \mbox{$6$ kn}. More details of the model parameters and the motion primitives used in the lattice-based planner are found in~\citep{bergman2020marine}.

The position-invariant cost function used in both steps of the framework is selected as 
\begin{equation}
\ell(\xbm, \ubm) = 1 + 0.1\nbm^{\sf T}\nbm + 100\dot{\nbm}^{\sf T}\dot{\nbm} + 100\dot{\alphabm}^{\sf T}\dot{\alphabm},
\vspace*{-0.5em}
\end{equation}
to trade-off between time, energy and smoothness. As for the cost related to the distance to obstacles, the value of the weighting parameter is selected as \mbox{$k_d = 1.5\cdot 10^{-3}$ } in \eqref{eq:lat-cost} and \eqref{eq:imp-cost}. Furthermore, the safety distance parameter is set to $d_{\mathrm{safe}} = 20 $ m which is approximately equal to the width of the ship's body.  
Finally, the COLREGs-related cost function is implemented as hard constraints, which corresponds to an infinite cost in \eqref{eq:lat-cost} and \eqref{eq:imp-cost}. The regions for each COLREGs situation are represented as polygons and are illustrated in Fig.~\ref{fig:colregs_borders}\footnote{The current implementation uses a fixed size of the polygons, but can be extended to depend on the velocity of the dynamic obstacle as in~\citep{tam2010collision}.}. The threshold parameters used to switch between COLREGs states in \eqref{eq:switch}-\eqref{eq:switch-back} are
\begin{equation}
\begin{aligned}
&\ubar{d} = 200, \; \; \ubar{t}_{\text{low}} = -10, \; \; \ubar{t}_{\text{high}} = 150, \\
&\bar{d} = 240, \; \; \bar{t}_{\text{low}} = -25, \; \; \bar{t}_{\text{high}} = 200.
\end{aligned}
\end{equation}

\begin{figure}[]
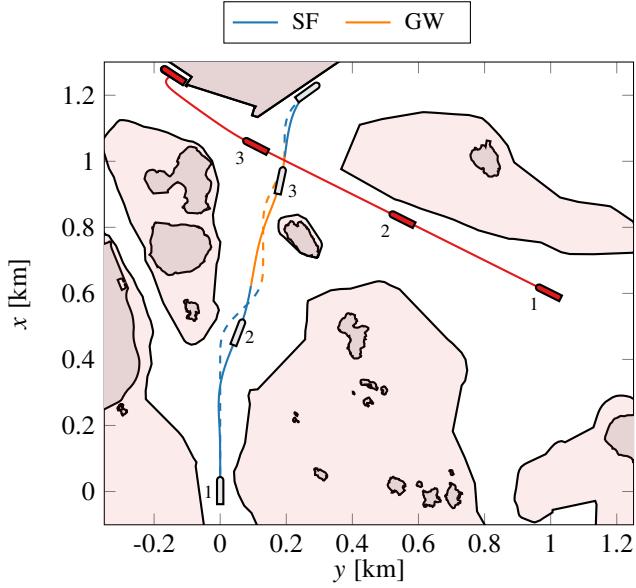

	\centering
	\setlength\figureheight{0.35\textwidth}
	\setlength\figurewidth{0.4\textwidth}
	\hspace*{2em}\input{label_giveway.tex} \\[0.2cm]
	\input{sim_giveway.tex} 
	\caption{ \label{fig:gw} \small (P3) A crossing scenario where the dynamic obstacle is approaching from the starboard side. The vessels are plotted at given timestamps (1-3) for interpretabilty. Path segments in orange represent parts where the ship is in state GW. } 
\end{figure}

\begin{figure}[t]
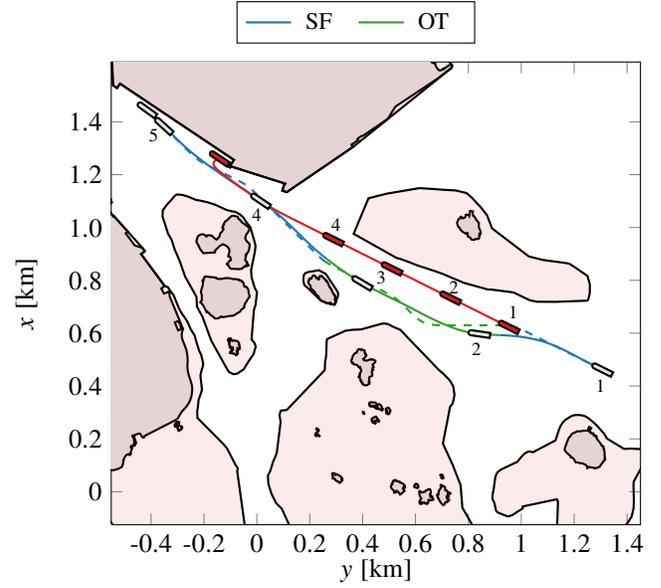

	\centering
	\setlength\figureheight{0.35\textwidth}
	\setlength\figurewidth{0.4\textwidth}
	\hspace*{2em}\input{label_overtaking.tex} \\[0.2cm]
	\input{sim_overtaking.tex} 
	\caption{ \label{fig:ot} \small (P4)  A scenario where the ship is able to overtake the dynamic obstacle. The vessels are plotted at given timestamps (1-5) for interpretabilty. Path segments in green represent parts where the ship is in state OT. } 
\end{figure}

The motion-planning scenarios are illustrated in Fig.~\ref{fig:sim_avoid}--\ref{fig:ot}, where four different trajectories are computed to different docking states with zero velocity. Since the performance of the proposed motion planner is evaluated in this section, the simulations are performed without any disturbances and the open-loop control law is used for evaluation. Moreover, the ship is aware of the future dynamic obstacle trajectories up to 150 seconds. After that, constant-velocity predictions are used for another 650 seconds.  

The result in Fig.~\ref{subfig-1:dummy} shows that the proposed motion planner is able to compute trajectories through narrow passages. In Fig.~\ref{subfig-2:dummy}, it can be seen that the motion planner is able to compute a trajectory that completely avoids a difficult COLREGs situation when another vessel is occupying the passage, by passing on the other side of the island. This is in contrast to methods that follows an initial reference computed without considering dynamic obstacles, such as~\citep{eriksen2020hybrid}, which instead react after the COLREGs situation has already occurred and hence is forced to wait until the other vessel has~passed. 

\begin{figure}[t!]
	\centering
	\hspace*{1.5em}\subfloat[$u$ (m/s) \label{subfig-1:d}]{%
		\setlength\figureheight{0.1\textwidth}
		\setlength\figurewidth{0.4\textwidth}
		\input{givewway_u.tex}\vspace*{-2em}
	} \\ \vspace*{-0.5em}
	\subfloat[$v$ (m/s)\label{subfig-2:d}]{%
		\setlength\figureheight{0.1\textwidth}
		\setlength\figurewidth{0.4\textwidth}
		\input{givewway_v.tex} 
	}\\ \vspace*{-1.5em}
	\hspace*{0.7em}\subfloat[$r$ (m/s)\label{subfig-3:d}]{%
		\setlength\figureheight{0.1\textwidth}
		\setlength\figurewidth{0.4\textwidth}
		\input{givewway_r.tex} 
	}
	\caption{ \small \label{fig:uvr-giveway} Linear velocity trajectories for the ship during the crossing scenario in Fig.~\ref{fig:gw}. Trajectory segments in orange represent parts where the ship is in state GW.    }
\end{figure}
In Fig.~\ref{fig:gw}, a crossing scenario (P3) is illustrated where the dynamic obstacle is approaching from the starboard side. This is identified correctly by the motion planner which switches to state GW and yields to the dynamic obstacle. It can be seen that the trajectory from the lattice-based motion planner is oscillatory. This behavior is explained by the following observations: 
\begin{enumerate}
	\item The small island prevents the planner from performing a course change to pass behind the dynamic obstacle.
	\item The lattice-based planner is limited to select motion primitives with velocities in $\mathcalOld{X}_d$, making it more beneficial to oscillate than to slow~down.
\end{enumerate}
The oscillatory behavior is then reduced by the improvement step, which can use all feasible velocities and is thus able to slow down instead of oscillating as seen in Fig.~\ref{fig:uvr-giveway}.

The last scenario (P4) is shown in Fig.~\ref{fig:ot}, where the motion planner finds a trajectory that safely overtakes the dynamic obstacle which is moving at a speed of 3 kn. A summary of the results from all scenarios is given in Table~\ref{tab:res_summary}. One notable result is the increase in computation time for the lattice-based planner for problem P2. An explanation to this increase is that the heuristic estimate in the HLUT is far from correct since it is unaware of the dynamic obstacle blocking the passage. Another result is that the improvement step is able to significantly reduce the aggregated applied force by the thrusters in all scenarios, and thus provide more smooth and energy-efficient trajectories compared to only using the lattice-based planner.  

\begin{table}[]
	\caption{\small Results for the motion planning problems in Fig.~\ref{fig:sim_avoid}--\ref{fig:ot}. The variable $t_{\text{lat}}$ represents the time  for the lattice-based planner to compute the initial trajectory, and $\tfinal$ is the time it takes to execute the trajectory. Finally, $F_{f}$ is the aggregated applied force for the initial lattice-based trajectory (lat) and after the improvement step (imp).    } \label{tab:res_summary}	
	\normalsize
	\centering 
	\setlength{\tabcolsep}{4pt}
	\begin{tabular}{ccccc}	
		Problem & $t_{\mathrm{lat}}$ [s] & $\tfinal$ [s] &$F_{f, \mathrm{lat}}$ [MN] & $F_{f, \mathrm{imp}}$ [MN]   \\
		\hline
		P1   & 0.11 & 533  & 71.3  & 52.7   \\
		P2   & 10.2 & 645  & 93.4 & 67.3 \\
		P3   & 0.24 & 467  & 67.8 & 46.4  \\ 
		P4   & 0.18 & 711  & 97.5 & 69.3  \\ 
		\hline
	\end{tabular}
\end{table}

\section{Conclusions and Future Work}\label{sec:conclusions}
This work proposes a two-step optimization-based motion planner for marine vessels that complies with the international regulations for preventing collisions at sea (\mbox{COLREGs}). It extends the work presented in~\citep{bergman2020marine} such that the motion planner is able to account for dynamic obstacles and the rules relevant for trajectory planning in COLREGs. In the first step, a lattice-based motion planner computes a suboptimal trajectory based on a finite, precomputed library of motion primitives. This step is here augmented with discrete states that represent the COLREGs situations of the ship with respect to other vessels. In the second step, direct optimal control techniques are applied to locally improve the trajectory from the lattice-based planner while keeping the COLREGs states fixed. The proposed motion planner is evaluated in a model of the Helsinki archipelago, where it successfully computes trajectories that complies with \mbox{COLREGs}. The results also show the motion planner's capability to completely avoid complex collision situations from occurring at all.

Future work includes to develop a heuristic function in the lattice-based planner that also captures both static and dynamic obstacles, to develop and evaluate the robustness of the implementation when the future behavior of the dynamic obstacles is uncertain, and to implement an emergency strategy.

\Urlmuskip=0mu plus 1mu
\bibliographystyle{abbrvnat}
{\small
	\bibliography{myrefs.bib}}		
\end{document}